\documentclass[a4paper,10pt]{article}

\usepackage{amsfonts,amssymb,mathrsfs,amscd}

\def \Z {{\mathbb {Z}}}

\def\eps{\varepsilon}

\usepackage[T1]{fontenc} 
\usepackage[cp1251]{inputenc}  
\usepackage[russian]{babel}

\title{\bf Slow convergence of Birkhoff ergodic averages}
\author{\bf Valery V. Ryzhikov}
\date{}
\begin{document}

\maketitle
\Large

\section{ Slow convergence of averages for $\Z$-actions}
Let $T$ be an ergodic automorphism of the probability space $(X,m)$, $f\in L_1(X,m)$. Birkhoff's theorem states that for almost all $x\in X$,
$$A(x,N,f):=\frac 1 n \sum_{i=1}^{N} f(T^ix) \to \int_X f\, dm,\ \ N\to\infty .$$
The rate of convergence depends on the automorphism and the function, which has given rise to many different problems, see \cite{P}.

The general result was obtained by Krengel \cite{K}: for an ergodic automorphism $T$, there is an arbitrarily slow rate of convergence of the Birkhoff averages, which is realized by choosing suitable functions $f\in L_1(X,m)$.
In the proposed note, we, modifying the approaches from \cite{23}, \cite{25}, realize the slow convergence of the Birkhoff averages for an arbitrary ergodic $\Z$-action. For this, we use the Rokhlin-Halmos lemma. The method is transferred to the actions of groups $\Z^n$ without any changes.
Its essence is as follows. Let the Birkhoff averages on time intervals of length $N_1$ for an ergodic automorphism
and a positive function $f$ for most $x$ become close to $\int_X f\, dm$. We set the function to zero on the union of trajectory segments of length $h\gg N_1$. Then outside these parts of the trajectories the Birkhoff averages will remain basically unchanged, but the integral of the function will change significantly. Thus, we have slowed down the rate of convergence at time $N_1$. In fact, the effect of the slowdown lasts for a very long time,
it depends on how large we have chosen $h$ compared to $N_1$.
By performing this procedure in a suitable way for a very rapidly growing sequence $N_k\to\infty$, we obtain the desired slowdown in the rate of convergence of the time Birkhoff averages for the final function.

Now let us present this more formally. Let $\eps_1>\delta_1>0$, a non-negative function $f_0\in L_1(X,m)$, $\|f_0\|=2$, a number $\eps>0$, and an ergodic automorphism $T$ of the probability space $(X,m)$ be given. By Birkhoff's theorem, there exists $N_1$
$$m\left(x:\ \left|\ A(x,N_1,f_0 - \int f_0dm \ \right|\ <\ \eps/10 \right)\ >1-\delta_1.$$
The Rokhlin-Halmos lemma implies the existence of a tower
$E_1= \bigsqcup_{i=1}^{N_1} T^iB_1$ such that $m(E_1)=\eps_1$.

If $N_1$ is large enough, then the set $C_1=X\setminus E_1$ satisfies the inequality
$$\left|\ \int f_0{\bf 1}_{C_1}dm - m(C_1)\int f_0dm \ \right|<\delta_1.$$
The inequality follows, for example, from the von Neumann ergodic theorem, taking into account the almost invariance of the set $C_1$. Indeed, the function $A(\cdot,M,f)$ is close to a constant for large $M$, therefore, due to the choice of $N_1\gg M$, the difference in the quantities
$$\int f_0{\bf 1}_{C_1}dm, \ \
\int A(\cdot,M,f){\bf 1}_{C_1}dm, \ \ m(C_1)\int f_0dm$$
can be made arbitrarily small (further we use $approx$).
Denoting $f_1=f_0{\bf 1}_{C_1}$, we note that
$$\int f_0dm-\int f_1dm \ >\ 0.9 \eps_1\int f_0dm.$$
For positive $\delta_1< \eps_1$ and arbitrarily large $N_1$, we obtain that for all $x$ outside the set of measure less than $\delta_1$, we have: on $C_1$,
$$A(x,N_1,f_1) = A(x,N_1,f_0)\approx \int f_0dm,$$ is satisfied, and on $E_1$ outside the mentioned set of measure less than $\delta_1$,
$$A(x,N_1,f_1) = 0$$
 is satisfied.
Thus, setting the function $f_0$ to be 0 on the set $E_1$  of measure $\eps_1$, we can realize for the new function $f_1$
deviations of $\eps/2$ (for most $x$) at a abitrary large time $N_1$:
$$m\left(x:\ \left|\ A(x,N_1,f_1) - \int f_1 dm\ \right|>\eps_1/2 \right)\ >1-\delta_1.$$
(By setting the function to zero on very long pieces of trajectories, we do not change most of the averages for $x$ outside these pieces at all, but the average of the function has changed significantly. For most $x$ in these pieces, the average has been set to zero, which gives the deviation we need on these $x$.)

Next, we choose an arbitrarily large $N_2$ so that
$$m\left(x:\ \left|\ A(x,N_2,f_1) - \int f_1dm \ \right|\ <\ \eps_2/10 \right)\ >1-\delta_2.$$
We consider a tower of the form $E_2=\bigcup_{i=1}^{h_2} T^iB_2$, $m(E_2)=\eps_2$. Let $$C_2=X\setminus (E_1\cup E_2), \ \ f_2=f_0{\bf 1}_{C_2}.$$
We choose the height $h_2\gg N_2$ for the tower $E_2$ so that $$m\left(x:\ \left|\ A(x,N_2,f_1) - \int f_2 dm\ \right|>\eps_2/2 \right)\ >1-\delta_2.$$

We put $C_{k}=X\setminus (E_1\cup\dots\cup E_{k})$ and $f_{k}=f_0{\bf 1}_{C_k}$. We choose the height $h_k\gg N_k$ so that
$$m\left(x:\ \left|A(x,N_k,f_k - \int f_k dm \right|>\eps_k/2 \right)\ >1-\delta_k.$$

For $$C=X\setminus \bigcup_{k= 1}^{\infty} E_{k}),
f=f_0{\bf 1}_C$$ we have
$$m\left(x:\ \left|A(x,N_k,f - \int fdm\right|>\eps_k/2 \right)\ >
1-2\sum_{i=k}^\infty \delta_i.$$

\vspace{2mm}
\bf Theorem 1. \it Let $T$ be an ergodic automorphism of the probability space $(X,m)$ and $f_0\in L_1(X,m)$ be a non-negative function,
$\|f_0\|>0$ and $a_i>0$ be given such that $\sum_{i}^\infty a_i<\infty$. For
any sequence $M_k\to +\infty$ there exist $N_k>M_k$ and a set $C$ of measure arbitrarily close to 1 such that for the function
$ f=f_0\,{\bf 1}_C$ we have
$$m\left(\ x:\ \left|\, A(x,N_k,f) - \int f\,dm \right|\ >\ a_k \ \right)\ \to\ 1.$$
Thus, for the function $f$ we have almost surely a given slow rate of convergence of Birkhoff averages to its integral. \rm

\section{ Slow convergence of averages for $\Z^n$-actions}
The method described above can be transferred without any fundamental changes to the case of actions of groups $\Z^n$. For simplicity, consider averaging over time (square) sets
$$Q_N=\{(z_1,\dots,z_n)\}\, :\, 1\leq z_1,\dots,z_n\leq N\}.$$

\vspace{2mm}
\bf Theorem 2. \it Let $\{T^z\}$ be an ergodic action of the group $\Z^n$ by automorphisms of the probability space $(X,m)$, $\sum_{i}^\infty a_i<\infty$, $a_i>0$. For
any sequence $M_k\to +\infty$ there exist $N_k>M_k$ and a function
$ f\in L_1(X,m)$ such that
$$m\left(\ x:\ \left|\, A(x,N_k,f) - \int f \, dm\, \right|\ >\ a_k \ \right)\ \to\ 1,$$ where
$$A(x,N,f) =\frac 1 {N^n} \sum_{z\in Q_N} f(T^{z}x).$$ \rm

\vspace{2mm}
In the proof of this theorem, Rokhlin's $\Z^n$-lemma is applied similarly to the case of $\Z$-actions. The role of towers on which we set the initial function to zero is now played by the sets $E_k= \bigsqcup_{z\in Q_{h_k}} T^zB_k.$
Instead of square averages, rectangular and more general configurations can be used. We leave the details as an exercise.
It is planned that a similar effect of slowing down the convergence rates for actions of countable amenable groups will be presented in a  note by I.V. Bychkov. For this generalization, an analogue of Rokhlin's lemma proposed by Ornstein and Weiss \cite{OW} could be used.
\section{ Медленная сходимость средних для $\Z$-действий}
Пусть  $T$ -- эргодический автоморфизм вероятностного пространства $(X,m)$,   $f\in L_1(X,m)$.  Теорема Биркгофа утверждает, что  для почти всех $x\in X$  выполнено  
$$A(x,N,f):=\frac 1 n \sum_{i=1}^{N} f(T^ix) \to  \int_X f\, dm,\ \ N\to\infty .$$ 
Скорость сходимости   зависит от автоморфизма и функции, в связи с чем  возникло множество разнообразных задач, см. \cite{P}.

Общий результат  был получен Кренгелем \cite{K}:  для  эргодического автоморфизма  $T$ имеет место    сколь угодно медленная скорость сходимости средних  Биркгофа, что  реализуется  при помощи выбора  подходящих функций  $f\in L_1(X,m)$. 
В предлагаемой заметке мы, модифицируя подходы  из работ \cite{23}, \cite{25},  реализуем  медленную  сходимость средних  Биркгофа для произвольного эргодического $\Z$-действия. Для этого  используем  лемму Рохлина-Халмоша. Метод без особых изменений переносится на действия групп $\Z^n$. 
Суть его в следующем. Пусть средние Биргофа на временных интервалах длины $N_1$  для эргодического автоморфизма 
и положительной функции $f$ для большинства $x$  стали  близкими к $\int_X f\, dm$.  Занулим функцию на объединении отрезков траекторий  длины $h\gg N_1$.
Тогда вне этих кусков траекторий средние Биркгофа в основном не изменятся, но существенно изменится  интеграл от функции. Тем самым мы замедлили скорость сходимости в момент времени $N_1$. На самом деле эффект замедления длится очень долго,
это зависит от того, насколько большим мы выбрали $h$  сравнению с $N_1$.  
Выполнив эту процедуру подходящим образом для очень быстро растущей последовательности $N_k\to\infty$,  получаем нужное замедление скорости сходимости временных средних Биркгофа для итоговой функции.

Теперь изложим сказанное более формально.  Пусть заданы $\eps_1>\delta_1>0$, неотрицательная  функция $f_0\in L_1(X,m)$,   $\|f_0\|=2$,  число  $\eps>0$ и эргодический автоморфизм $T$ вероятностного пространства $(X,m)$.    В силу теоремы Биркгофа найдется $N_1$ 
$$m\left(x:\ \left|\ A(x,N_1,f_0 - \int f_0dm \ \right|\ <\ \eps/10 \right)\ >1-\delta_1.$$ 
Из леммы Рохлина-Халмоша вытекает существование башни 
$E_1= \bigsqcup_{i=1}^{N_1} T^iB_1$  такой, что  $m(E_1)=\eps_1$. 
 
 Если $N_1$ достаточно велико, то  для множества $C_1=X\setminus E_1$  выполнено неравенство
$$\left|\ \int  f_0{\bf 1}_{C_1}dm - m(C_1)\int f_0dm \ \right|<\delta_1.$$
Неравенство следует, например,  из эргодической теоремы фон Неймана  с учетом   почти инвариантности множества $C_1$. Действительно, функция  $A(\cdot,M,f)$ близка к константе для большого $M$, поэтому благодаря выбору   $N_1\gg M$  различие величин
$$\int  f_0{\bf 1}_{C_1}dm, \ \ 
\int A(\cdot,M,f){\bf 1}_{C_1}dm, \ \ m(C_1)\int f_0dm$$ 
можно сделать произвольно малым (далее используем $approx$).
Обозначая $f_1=f_0{\bf 1}_{C_1}$, замечаем, что 
$$\int  f_0dm-\int  f_1dm \ >\ 0.9 \eps_1\int f_0dm.$$ 
Для положительного $\delta_1< \eps_1$ и сколь угодно большого $N_1$, получим, что для всех $x$ вне   множества меры, меньшей $\delta_1$, имеем: на $C_1$ выполнено 
$$A(x,N_1,f_1) =  A(x,N_1,f_0)\approx \int f_0dm,$$ а на $E_1$ вне упомянутого множества меры, меньшей $\delta_1$,
 выполнено  $$A(x,N_1,f_1) =  0.$$
Таким образом, занулив функцию $f_0$ на множестве $Е_1$ меры $\eps_1$, мы можем  реализовать для полученной функции $f_1$
отклонения  на $\eps/2$ (для большинства $x$)  в сколь угодно далекий момент времени $N_1$:
$$m\left(x:\ \left|\ A(x,N_1,f_1) - \int f_1 dm\ \right|>\eps_1/2  \right)\ >1-\delta_1.$$
Неформальное пояснение:  зануляя функцию на очень длинных кусках траекторий,
мы большинство средних для $x$  вне этих кусков вообще не меняем, но при этом существенно изменилось среднее функции.  Для большинство  $x$  в этих кусках среднее обнулилось, что дает нужное нам отклонение на этих $x$.

Далее выбираем сколь угодно большое  $N_2$ так, чтобы 
$$m\left(x:\ \left|\ A(x,N_2,f_1) - \int f_1dm \ \right|\ <\ \eps_2/10 \right)\ >1-\delta_2.$$ 
Рассматриваем башню вида  $E_2=\bigcup_{i=1}^{h_2} T^iB_2$, $m(E_2)=\eps_2$.  Положим $$C_2=X\setminus (E_1\cup E_2), \ \ f_2=f_0{\bf 1}_{C_2}.$$ 
Высоту  $h_2\gg N_2$ для башни $E_2$  выбираем так, чтобы $$m\left(x:\ \left|\ A(x,N_2,f_1) - \int f_2 dm\ \right|>\eps_2/2  \right)\ >1-\delta_2.$$  

Положим $C_{k}=X\setminus (E_1\cup\dots\cup E_{k})$ и   $f_{k}=f_0{\bf 1}_{C_k}$. Высоту  $h_k\gg N_k$ выбираем так, чтобы 
$$m\left(x:\ \left|A(x,N_k,f_k - \int f_k dm \right|>\eps_k/2  \right)\ >1-\delta_k.$$  

Для  $$C=X\setminus \bigcup_{k= 1}^{\infty} E_{k}), 
f=f_0{\bf 1}_C$$ имеем 
$$m\left(x:\ \left|A(x,N_k,f - \int fdm\right|>\eps_k/2  \right)\ >
1-2\sum_{i=k}^\infty \delta_i.$$  

\vspace{2mm} 
\bf Теорема 1. \it Пусть  $T$ -- эргодический автоморфизм вероятностного пространства $(X,m)$ и   $f_0\in L_1(X,m)$ -- неотрицательная  функция,
$\|f_0\|>0$ и заданы  $a_i>0$ такие, что $\sum_{i}^\infty a_i<\infty$.  Для 
всякой последовательности $M_k\to +\infty$ найдутся  $N_k>M_k$ и множество  $С$ меры, сколь угодно близкой к 1, такое, что  для функции 
$ f=f_0\,{\bf 1}_C$  выполнено 
$$m\left(\ x:\ \left|\, A(x,N_k,f) - \int f\,dm \right|\ >\ a_k \ \right)\ \to\ 1.$$ 
Таким образом, для функции $f$ мы  имеем  почти наверное заданную  медленную скорость сходимости средних Биркгофа к ее интегралу. \rm 

\section{  Медленная сходимость средних для $\Z^n$-действий}
Метод, изложенный выше,  без принципиальных изменений  переносится на случай действий групп $\Z^n$. Рассмотрим для простоты усреднения вдоль временных  (квадратных)  множеств 
$$Q_N=\{(z_1,\dots,z_n)\}\, :\, 1\leq z_1,\dots,z_n\leq N\}.$$

\vspace{2mm} 
\bf Теорема 2. \it Пусть  $\{T^z\}$ -- эргодическое действие группы $\Z^n$ автоморфизмами  вероятностного пространства $(X,m)$,   $\sum_{i}^\infty a_i<\infty$,  $a_i>0$.  Для 
всякой последовательности $M_k\to +\infty$ найдутся  $N_k>M_k$ и   функция 
$ f\in L_1(X,m)$  такие, что выполнено 
$$m\left(\ x:\ \left|\, A(x,N_k,f) - \int f \, dm\, \right|\ >\ a_k \ \right)\ \to\ 1,$$  где 
$$A(x,N,f) =\frac 1 {N^n} \sum_{z\in Q_N} f(T^{z}x).$$  \rm

\vspace{2mm} 
В доказательстве этой теоремы  аналогично случаю $\Z$-действий применяется $\Z^n$-лемма Рохлина. Роль башен, на которых мы зануляем иходную функцию,  теперь играют множества $E_k= \bigsqcup_{z\in Q_{h_k}} T^zB_k.$
Вместо квадратных усреднений можно использовать  прямоугольные и более общие конфигурации. Оставляем детали в качестве упражнения. 
 Планируется, что  подобный эффект замедления скоростей сходимости для  действий счетных аменабельных групп будет изложен в  отдельной  заметке И.В. Бычкова. Для этого обобщения  будет использован  аналог леммы Рохлина, предложенный Орнстейном и Вейсом \cite{OW}.

\end{document}